\documentclass[12pt]{amsart}

\newtheorem{theorem}{Theorem}[section]
\newtheorem{proposition}[theorem]{Proposition}
\newtheorem{lemma}[theorem]{Lemma}

\theoremstyle{definition}
\newtheorem{remark}[theorem]{Remark}

\usepackage{graphics}
\usepackage{t1enc}
\usepackage{textcomp}
\usepackage{amsmath, amstext,amscd, amsthm}
\usepackage{epic} \usepackage{eepic}
\usepackage{amssymb}

\def\Spec{{\rm Spec}\,}

\def\tens{\otimes}

\def\Frac#1#2{{\displaystyle{{#1} \overwithdelims.. {#2}}}}

\def\wt{\widetilde{\tau}}
\def\wX{\widetilde{X}}
\def\wJ{\widetilde{J}}

\def\Z{\mathbb{Z}}

\def\C{\mathbb{C}}
\def\P{ \mbox{P\hspace{-.3em}I}}

\def\O{\mathcal{O} }
\def\P{\mathbb{P}}
\def\Supp{\text{Supp}\,}

\makeatletter

\@addtoreset{equation}{section} \makeatother

\begin{document}

\baselineskip=15.5pt

\title[Frobenius action on vector bundles
in small odd characteristics]{The Frobenius action on rank $2$
vector bundles over curves in small genus and small characteristic}

\author[L. Ducrohet]{Laurent Ducrohet}

\address{CMLS, Ecole Polytechnique, 91128 Palaiseau Cedex, France}

\email{ducrohet@math.polytechnique.fr}

\date{}

\maketitle

\begin{abstract}Let $X$ be a general proper and smooth curve of genus 2 (resp. of
genus 3) defined over an algebraically closed
field of characteristic $p$. When $3\leq p \leq 7$, the action of
Frobenius on rank 2 semi-stable vector bundles with trivial
determinant is completely determined by its restrictions to the 30
lines (resp. the 126 Kummer surfaces) that are invariant under the
action of some order 2 line bundle over $X$. Those lines (resp.
those Kummer surfaces) are closely related to the elliptic curves
(resp. the abelian surfaces) that appear as the Prym
varieties associated to double \'etale coverings of $X$. We are therefore able to compute the explicit equations defining Frobenius action in
these cases. We perform some of these computations and draw some geometric
consequences.
\end{abstract}

\section{Introduction}

Let $k$ be a algebraically closed field of positive characteristic
$p$ and let $X$ be  an irreducible,  proper and smooth curve of
genus $g$ over $k$. Let $X_s$ ($s\in \Z$ since $k$ is perfect) be
the $p^s$-twist of $X$ and let $J$ (resp. $J_s$) denote its Jacobian
variety (resp. the $p^s$-twist of its Jacobian variety). Also, let
$\Theta$ denote a symmetric principal polarization for $J$ (associated to a theta characteristic $\kappa_0$). Denote
by $S_X(r)$ the (coarse) moduli space of semi-stable rank $r$ vector
bundles with trivial determinant over $X$. The map $E\mapsto F_{\rm
abs}^*\, E$ defines a rational map $S_{X}(r) \dashrightarrow
S_{X}(r)$ the $k$-linear part $V_r : S_{X_1}(r) \dashrightarrow
S_{X}(r)$ of which is called the (generalized) Verschiebung.

Our interest in this situation stems from the fact (see \cite{LS})
that a stable rank $r$ vector bundle $E$ over $X$ corresponds to an
(irreducible) continuous representation of the algebraic fundamental
group $\pi_1(X)$ in ${\rm GL}_r(\bar{k})$ (endowed with the discrete
topology) if and only if one can find an integer $n>0$ such that
${F_{\rm abs}^{(n)}}^*E\cong E$. Thus, natural questions about the
generalized Verschiebung $V_r : S_{X_1}(r) \dashrightarrow S_{X}(r)$
arise like, e.g., its surjectivity, its degree, the density of
Frobenius-stable bundles (i.e., those vector bundles whose pull-back
by Frobenius iterates are all semi-stable), the loci of
Frobenius-destabilized bundles...

For general $(g,\,r,\,p)$, not much seems to be known (see the
introductions of \cite{LP1} and \cite{LP2} for an overview of this subject).\\

We will focus on the rank $r=2$ in low genus ($g=2$ or 3) and we
will let $S_X$ stand for $S_X(2)$. When $k=\C$, Narasimhan and
Ramanan have given explicit descriptions for $S_X$ as a subvariety
of $|2\Theta|$. Namely, $S_X$ is isomorphic to $|2\Theta|$ in the
genus 2 case (see \cite{NR1})  and $S_X$ identifies with the Coble
quartic surface in $|2\Theta|$ in the genus 3 non-hyperellitic case
(see \cite{NR2}). These descriptions also hold over any
algebraically closed field of odd characteristic and $V_2$ lifts to
a rational map $\widetilde{V} : |2\Theta_1| \dashrightarrow
|2\Theta|$ (see Proposition \ref{widetildeVexists} for the genus 3
case) given by a system of $J[2]$-equivariant homogeneous
polynomials of degree $p$.

In genus 2 and characteristic $p=3$, Laszlo and Pauly gave in
\cite{LP2} the cubic equations of $\widetilde{V}$ in showing that
this rational map coincided with the polar map of a quartic surface,
isomorphic to ${\rm Kum}_X$, embedded in $|2\Theta_1|$. The proof
uses the fact that the action of Frobenius is equivariant under the
action of $J[2]$ in odd characteristic as well as a striking
relationship (see \cite{vG}) between cubics and quartics on
$|2\Theta_1|$ that are invariant under the action of $J[2]$. For a
general odd $p$, the base locus of $\widetilde{V}$ coincides at
least set-theoretically with the locus of Frobenius destabilized
bundles (i.e., those stable vector  bundle $E$ such that $F^*E$ is
unstable) and it has been much studied (see \cite{LnP} and
\cite{Os1}).

In this article, we shall suppose that the characteristic $p$ of the
base field is odd. Given a line bundle $\tau$ of order 2 over $X$
and a $\tau$-invariant (i.e., satisfying $E\tens \tau \xrightarrow
\sim E$) semi-stable degree 0 vector bundle $E$, one can give $E$ a
structure of invertible $\O_X\oplus \tau$-module. In other words, if
$\pi : \widetilde{X} \to X$ is the degree 2 étale cover
corresponding to $\tau$, there is a degree 0 line bundle $L$ over
$\widetilde{X}$ such that $E\cong \pi_*(L)$. Because $\pi$ is étale,
one has $F_{\rm abs}^*(\pi_*(L))\cong \pi_*(F_{\rm abs}^*(L))$.
Furthermore, requiring $E$ to have trivial determinant forces $L$ to
be in some translate of the Prym variety $P_\tau$ associated to
$\pi$ (which has genus $g-1$). On the one hand, the associated
morphism $P_\tau \to |2\Theta|^\tau$ factors through the Kummer
morphism $P_\tau \twoheadrightarrow P_\tau/\{\pm\}$ and, on the
other hand, as multiplication by $p$ over an abelian variety
commutes with the inversion, it induces an endomorphism of
$P_\tau/\{\pm\}$. If $g=2$, the Prym varieties are elliptic curves
and there are formulae (see \cite{Si}) that allow us to compute
explicitly the $k$-linear part $\widetilde{V}_\tau : \P^1 \to \P^1$
of this endomorphism. If $g=3$, the Prym varieties are Jacobian
varieties of genus 2 curves and one can use the genus 2 case
results.

Through representations of Heisenberg groups, we prove the key
result of this article
\begin{theorem} \label{MainTheorem} Let $k$ be an algebraically closed field of
characteristic $p=3,\,5$ or 7. Let $X$ be a  smooth, proper, general
curve of genus 2 or 3 over $k$. There is a rational map
$\widetilde{V} : |2\Theta_1| \dashrightarrow |2\Theta|$ extending
the generalized Verschiebung $V_2 : S_{X_1}\dashrightarrow S_X$ that
is completely determined by its restrictions $\widetilde{V}_\tau :
|2\Theta_1|^{\tau_1} \dashrightarrow |2\Theta|$ to the
$\tau_1$-invariant locus of $|2\Theta_1|$, $\tau_1$ ranging in the
non zero elements of $J_1[2]$.
\end{theorem}

Therefore, one can explicitly compute the equations of $V_2 :
S_{X_1} \dashrightarrow S_X$ and we perform these computations in
genus 2 and characteristic 3,\,5 and 7, as well as in genus 3 and
characteristic 3.

In the genus 2 case, we give the following generalization of the
results of \cite{LP2} in characteristic 3 :

\begin{proposition}\label{LPgeneralizedg2} Let $X$ be a general, proper and smooth curve
of genus 2 over an algebraically closed field of characteristic $p$.
For $p=3,\,5$ or $7$, there is a degree $2p-2$ irreducible
hypersurface $H$ in $|2\Theta_1|$ such that the equality of divisors
in $|2\Theta_1|$
$$\widetilde{V}^{-1}({\rm Kum}_{X})= {\rm Kum}_{X_1}+2 H $$ holds
\emph{scheme-theoretically}.
\end{proposition}

Randomly choosing curves, the following pattern arises, analogous to
the characteristic 3 case : The base locus of $V$ is strictly
contained in the singular locus of $H$. The latter has dimension 0,
is contained in the stable locus of $S_{X_1}$, as well as in the
inverse image of the singular points of ${\rm Kum}_X$ which is
1-dimensional. Unfortunately, the Groebner basis computation
required to check this statement for the generic curve seems too
heavy and we could not check this result globally.\\

In the genus 3 case, the results of \cite{LP2} in characteristic 3
generalize as follows :

\begin{theorem} \label{LPgeneralizedg3} Assume that $X$ is a general, smooth and projective curve
of genus 3 over an algebraically closed field of characteristic 3.\\
There is an embedding $\alpha : {\rm Cob}_{X} \hookrightarrow
|2\Theta_1|$ such that the cubic equations of the rational map
$\widetilde{V} : |2\Theta_1| \dashrightarrow |2\Theta|$ lifting the
generalized Verschiebung $V_2 : S_{X_1} \dashrightarrow S_X$ are
given by the 8 partial derivatives of the quartic equation of
$\alpha({\rm Cob}_X) \subseteq |2\Theta_1|$. In other words,
$\widetilde{V}$ is the polar map of the hypersurface $\alpha({\rm
Cob}_X)$.\\ In particular, the base locus (equivalently, the locus
of Frobenius destabilized bundle) is the intersection $\alpha({\rm
Kum}_X)\bigcap {\rm Cob}_{X_1}$.
\end{theorem}

All the computations have been carried out using MAGMA Computational
Algebra System, on the servers MEDICIS hosted at the Ecole
Polytechnique.

I would like to thank Y. Laszlo for having introduced me to this
question, for his help and encouragements.\\

\section{Vector bundles and Theta group representations}

\subsection{Action of $J[2]$ on the moduli space $S_X$}
Following \cite[1.8]{Ra}, there is a
morphism $D : S_X \to |2\Theta|$ mapping a (S-equivalence class
of) semi-stable rank 2 vector bundle with trivial determinant $E$ to
the unique effective divisor in $|2\Theta|$ with support the set
\begin{equation}\label{NRBmap}
\Supp \, D[E]= \{j\in J|\, H^0(X,\,E\tens j\tens \kappa_0)\neq
0\}\end{equation} We will consider the morphism $b : J \to S_X$ defined by
$j \mapsto [j\oplus j^{-1}]$,
the Kummer morphism ${\rm K}_X : J \to |2\Theta|$ which maps $J$ onto the
Kummer variety ${\rm Kum}_X\cong J/\{\pm 1\}$ and the
morphism $\varphi_{2\Theta} : J \to |2\Theta|^*$
associated to $\O(2\Theta)$. Also, following \cite[Sect. 2]{Be}, introduce the subvariety $\Delta$ of $S_X$ with support
the set
\begin{equation*}\{[E] \in S_X|\, H^0(E\tens \kappa_0)\neq 0\}\end{equation*}

\begin{proposition} \label{Beauville} \textup{(1)} One has $b^*(\O(\Delta)) \cong
\O(2\Theta)$ and the map $b^* : H^0(S_X,\,\O(\Delta)) \to
H^0(X,\,\O(2\Theta))$ is an isomorphism.\\
Identifying $|2\Theta|^*$ and $|\Delta|^*$ via this isomorphism, the
 morphism $\varphi_\Delta : S_X \to |2\Theta|^*$ associated
 to the linear system $|\Delta|$ gives a commutative diagram
\begin{equation}\label{beauvillediag}\unitlength=0.6cm
\begin{picture}(15,9)
\put(0,4){$J$} \put(5.6,4){$S_X$} \put(10,7.5){$|2\Theta|^*$}
\put(10,0.5){$|2\Theta|$}

\put(10.6,7.2){\vector(0,-1){6}}

\put(7.4,3.7){\vector(1,-1){2.5}}

\put(7.4,4.8){\vector(1,1){2.5}}

\put(0.8,4.2){\vector(1,0){4}}

\qbezier(0.8,4.5)(4.7,7.7)(9.7,7.7) \put(9.7,7.7){\vector(1,0){0.1}}

\qbezier(0.8,3.9)(4.7,0.7)(9.7,0.7) \put(9.7,0.7){\vector(1,0){0.1}}

 {\footnotesize
\put(10.7,4){$\wr$} \put(3.2,4.5){$b$} \put(5,1.6){${\rm K}_X$}
\put(5,6.4){$\varphi_{2\Theta}$}
 \put(8.9,2.7){$D$} \put(8.9,5.7){$\varphi_{\Delta}$ } }
\end{picture}\end{equation}
(where the vertical arrow is Wirtinger's isomorphism) \cite[Proposition 2.5]{Be}.\\
\textup{(2)} If $g=2$, the morphism $D$ is an isomorphism \cite{NR1}.\\
\textup{(3)} If $g=3$ and $X$ is not hyperelliptic, $D$ is a closed
immersion whose image is the Coble quartic surface ${\rm Cob}_X$ \cite{NR2}.\end{proposition}
Define actions of $J[2]$ on $S_X$ and $|2\Theta|$ (hence on $|2\Theta|^*$ by duality) respectively
by $(\tau,\,[E]) \mapsto [E \tens
\tau] $,
$(\tau,\,D) \mapsto T_\tau^*D$ (where
$T_\tau$ is the translation by $\tau$ on the Jacobian $J$)
All the maps in the
diagram above are $J[2]$-equivariant (see
\cite[Remark 2.6]{Be}).

\subsection{Theta groups and representations}\label{Thetagroups}
Let $A$ be any abelian variety over $k$ and $L$ an ample line bundle $A$. Following \cite[Sect. 1]{Mu2}, let $\mathcal{G}(L)$
(resp. $K(L)$) be the group scheme (resp. the finite group scheme)
such that, for any $k$-scheme $S$,
$$\mathcal{G}(L)(S)=\{(x,\,\gamma)|\,x\in A(S),\ \gamma : L
\xrightarrow \sim T_x^*L \}\ \ ({\rm resp.}\  K(L)(S)=\{x\in
A(S)|\,T_x^*L \cong L\})$$ The commutator in the theta group $\mathcal{G}(L)$ induces a non degenerate skew-symmetric
bilinear form denoted by $e^L : K(L) \times K(L) \to \mathbb{G}_m$.
Suppose that $K(L)$ is reduced-reduced, i.e., $K(L)$ is reduced and
its Cartier dual is also reduced, then $L$ is said to be of
separable type.

The ample line bundle
$\O(2\Theta)$ over $J$ is of that kind and we introduce the following notation :
\begin{equation*} W := H^0(J,\,\O(2\Theta)),\
\mathcal{G}(2) := \mathcal{G}(\O(2\Theta)),\
e_{2,\,J} : = e^{\O(2\Theta)} : J[2] \times J[2] \to \mu_2\subset \mathbb{G}_m\end{equation*}
The vector space
 $W$ is the unique (up to isomorphism)
irreducible representation of weight 1 of $\mathcal{G}(2)$ \cite[Sect. 1]{Mu2}. By
duality, there is an action (of weight -1) of $\mathcal{G}(2)$ on
$W^*$.

 We write $H=(\Z/2\Z)^g$ and $\hat{H}=\text{Hom}((\Z/2\Z)^g,\,k^*)$
 (we
identify $H$ and $\hat{H}$ by means of the bilinear  form
$(\alpha,\,\beta) \mapsto {^t\alpha.\beta}$ with values in
$\mathbb{F}_2$) and we consider the associated \emph{Heisenberg group} $\mathcal{H}$ with underlying set $k^* \times H \times \hat{H}$. We let $E$ denote
the non-degenerate bilinear form on $H\times \bar{H}$ defined by the commutator in
$\mathcal{H}$. Recall that a  \emph{theta structure} $\widetilde{\phi} :
\mathcal{H} \xrightarrow \sim \mathcal{G}(2)$ on
$\mathcal{G}(2)$ is entirely determined by the images of $H$ and $\bar{H}$ and a \emph{theta basis} $\{X_\alpha|\, \alpha\in H\}$ of $W$ is canonically (up to multiplicative scalar) associated to it. It satisfies the following properties :
$$\beta.X_\alpha=X_{\alpha+\beta} \text{ for any } \alpha,\,\beta \in H,\ \alpha^*.X_\alpha=\alpha^*(\alpha)X_\alpha \text{ for any } \alpha,\,\alpha^* \in
H\times \hat{H}.$$

\subsection{$\tau$-invariant vector bundles and \'etale double
covers} \label{tauinvlines-prymvars}
Choose a non-zero element $\tau$ of $J[2]$ and consider the associated double \'etale cover $\pi : \widetilde{X}:=\textbf{Spec}(\O_X\oplus
\tau) \to X$ with genus $2g-1$ (Hurwitz formula).
Letting $\wJ$ denote the Jacobian of $\wX$, denote by \begin{equation}\label{ellipticPrymPtau} P_\tau:=\ker ({\rm Nm} : \wJ \to J)^0\end{equation}
the Prym variety associated to $\pi$, defined as the neutral component of $\ker {\rm Nm}$ (see \cite{Mu3} for general properties of Prym varieties). The homomorphism
$\sigma :J \times P_\tau \to \widetilde{J}$
defined set-theoretically by $(j,\,L) \mapsto \sigma(j,\,L) :=
\pi^*(j)\tens L$ has reduced-reduced kernel $K_\sigma$. As $\pi$ is \'etale and
$\pi^*(z)^2\cong \O_{\wX}$, $\pi^*(\kappa_0\tens
z)$ is a theta characteristic for $\wX$. Denote by
$\widetilde{\Theta}_z$ the corresponding symmetric principal divisor on $\widetilde{J}$. One has the set-theoretical
equality
\begin{equation*}
{\rm Supp}\, D[\pi_*L\tens z]= (\pi^*)^{-1}{\rm
Supp}\,(T^*_L\widetilde{\Theta}_z)
\end{equation*}

\begin{proposition} \label{dtauz} \textup{(1)} Choose an element $z$ in $S_\tau=\{z\in
J|\,z^2=\tau\}$. Then, there is a well-defined morphism
\begin{equation*} d_{\tau,\,z} : P_\tau \to S_X \end{equation*}
mapping $L$ to $[\pi_*L \tens z]$. Furthermore, if $\pi_*L \tens z$
is
strictly semi-stable, then $L$ is in $P_\tau[2]$.\\
 A vector
bundle of this form is $\tau$-invariant, i.e., it is equipped with
an isomorphism $(\pi_*L \tens z) \tens \tau \xrightarrow \sim \pi_*L
\tens z$.\\
Conversely, any $\tau$-invariant semi-stable rank 2 vector bundle
with trivial determinant is S-équivalent to $\pi_*L \tens z'$ for
some $L$ in $P_\tau$ and some $z'$ in $S_\tau$.\\
\textup{(2)} There is a well-defined morphism
\begin{equation*}\delta_{\tau,\,z} : P_\tau \to |2\Theta|\end{equation*}
mapping $L$ to the divisor
$(\pi^*)^{-1}(T^*_L\widetilde{\Theta}_z)$.\\
\textup{(3)} The morphism $\delta_{\tau,\,z}$ agrees with the composite $D\circ d_{\tau,\,z}$.
\end{proposition}
\begin{proof} (1) Using projection formula,
one finds that the degree zero rank 2 vector bundle $\pi_*L$ is semi-stable and non-stable if and only if $(M,\,L^{-1})$ lie in $K_\sigma$, i.e., if
$L\cong \pi^*M$ is in $P_\tau[2]$. Those statements hold
after further imposing that $L$ lies in $P_\tau$ and tensoring by
$z$. In this case, $\det (\pi_*L\tens z)\cong {\rm Nm}(L)\tens \tau \tens z^2$
is trivial. The Poincar\'e line
bundle over $\wX \times \wJ$ provides a
family of rank 2 and trivial determinant vector bundles over $X$ parameterized by $P_\tau$ and the coarse moduli property induces the morphism
$d_{\tau,\,z}$. Because $\pi^*\tau $ is trivial, the projection
formula ensures that $\pi_*L\tens z$ is $\tau$-invariant.
Conversely, a $\tau$-invariant vector bundle with rank 2 and trivial determinant is isomorphic to $\pi_*L$ for some line bundle $L$ on $\wX$ and the statement follows for determinant reasons.\\
(2) follows from \cite{Mu3} and (3) is clear since a divisor in $|2\Theta|$ is entirely determined by its support.\end{proof}

\subsection{Theta groups and Prym varieties} The Prym varieties $P_\tau$ is actually a principally polarized
abelian variety (see \cite{Mu3}, Sections 2 and 3). Choosing any symmetric principal divisor $\Xi$ on $P_\tau$, the line bundle $\O(2\Xi)$ is canonical
and one has
$$\sigma^*(\O(\widetilde{\Theta}_z))\cong \O(2\Theta) \boxtimes
\O(2\Xi)$$ Therefore, \cite[Sect. 23, Thm. 2]{Mu1} ensures that there is a unique level subgroup $K_\sigma \cong \widetilde{K}_\sigma $ in the Heisenberg group $
\mathcal{G}(\O(2\Theta) \boxtimes \O(2\Xi))$ such that $\sigma_*(\O(2\Theta) \boxtimes
\O(2\Xi))^{\widetilde{K}_\sigma} \cong \O(\widetilde{\Theta}_z)$.
Denote by $\wt:=\wt (z)$ the image of $\tau$ via the lifting
$K_\sigma \xrightarrow{\sim} \widetilde{K}_\sigma $. It follows from \cite[Sections 4 and 5]{Mu3} that there is a unique (up to multiplicative constant) isomorphism
\begin{equation}\label{isomchi} \chi :
H^0(J,\,\O(2\Theta))^{<\wt>} \xrightarrow{\sim}
H^0(P_\tau,\,\O(2\Xi))^*\end{equation} of Heisenberg
representations and that the morphism $\delta_{\tau,\,z} :
P_\tau \to |2\Theta|$
factors as the composite
$$ P_\tau \xrightarrow{\varphi_{2\Xi}}\P
H^0(P_\tau,\,\O(2\Xi))^* \xrightarrow{\sim} \P
H^0(J,\,\O(2\Theta))^{\wt} \subset |2\Theta|$$ where the (canonical)
isomorphism is deduced from $\chi$.

\begin{remark}\label{changingz} Notice that the set $S_\tau$
is principal homogeneous under the action $J[2]$ and check that  $\wt(z+\alpha)=e_{2,\,J}(\alpha,\,\tau)\wt(z)$. In particular,
$\wt(z)=\wt(-z)$.
\end{remark}

\section{The action of Frobenius}

\subsection{Theta groups in characteristic $p$} For any scheme $S$ over $k$, we
introduce the $p$-twist $S_1$ of $S$ and the (relative) Frobenius $F
: S \to S_1$ (which is $k$-linear by contrast with $F_{\rm abs}$), both defined by the following commutative diagram
\begin{center}
\unitlength=0.6cm
\begin{picture}(12,6)
\put(0,4.5){$S$}

 \dashline{0.1}(0.8,4.3)(2.6,3.1) \put(2.6,3.1){\vector(3,-2){0.1}}

\qbezier(1,4.6)(5.2,4.6)(7.9,3) \put(7.9,3){\vector(1,-1){0.1}}

\qbezier(0.2,4.2)(0.2,1.2)(2.2,0.6)
\put(2.2,0.6){\vector(1,-1){0.1}}

\put(3.3,2.5){$S_1$}

\put(2.4,0){ \put(6,2.5){$S$} \put(0.4,0){$\Spec k
$}\put(5.4,0){$\Spec k$} \put(2.3,2.7){\vector(1,0){3.1}}
\put(2.5,0.2){\vector(1,0){2.5}} \put(1.2,2.3){\vector(0,-1){1.5}}
\put(6.2,2.3){\vector(0,-1){1.5}} \small \put(3.2,0.4){$F_{\rm
abs}$} \put(3.5,2.9){$i$}  \put(3.5,1.2){$\square$} }

\small \put(2,3.7){$F$} \put(5.1,4.3){$F_{\rm abs}$}
\end{picture}\end{center}
where the square is cartesian. Apart from $J_1$ and $S_X(r)_1$, define $\Theta_1$ (resp. $\Delta_1$) as the $p$-twists
$i^*\Theta$ (resp. $i^*\Delta$). As before, define
\begin{equation*} W_1 := H^0(J_1,\,\O(2\Theta_1)),\
\mathcal{G}_1(2) := \mathcal{G}(\O(2\Theta_1)),\
e_{2,\,J_1} : = e^{\O(2\Theta_1)}\end{equation*}

Assume from now on that $X$ is an ordinary curve. The line bundle  $\O(p\,\Theta_1)$ is no
longer of separable type but Sekiguchi proved in \cite{Se} that the
main results about theta groups proved by Mumford in case of line bundles of separable type can be extended in that
case (we refer to \cite[Section 2]{LP1}). Because $p$ is odd, $[p\,]$ induces identity on $J[2]$ and the
restrictions $F : J[2] \to J_1[2]$ and $V : J_1[2] \to J[2]$ are
isomorphisms, inverse one to each other, whence
 a natural
action of $J[2]$ on the spaces $S_{X_1}$, $|2\Theta_1|$...

\subsection{Finding equations for the Frobenius action on vector
bundles} In the case of vector bundles, denote by $V_r : S_{X_1}(r)  \dashrightarrow
S_X(r)$ and call
(generalized) Verschiebung
the pull-back by $F :X \to X_1$ .
The diagram
\begin{equation}\label{VS_XVJcommute}
\unitlength=0.6cm
\begin{picture}(10,4)
\put(0.6,0){$S_{X_1}$} \put(8.6,0){$S_{X}$} \put(0.8,3){$J_1$}
\put(8.9,3){$J$} \put(2,3.2){\vector(1,0){6}}
\put(1.2,2.7){\vector(0,-1){1.8}} \put(9,2.7){\vector(0,-1){1.8}}
\put(7.6,0.2){\vector(1,0){0.1}} \dashline{0.2}(2.6,0.2)(7.6,0.2)

\small \put(5,3.4){$V$} \put(5,0.6){$V_2$} \put(0.4,1.5){$b_1$}
\put(9.3,1.5){$b$}
\end{picture}
\end{equation}
commutes and that means that  $V_2 : S_{X_1}
\to S_X$ extends $V : J_1 \to J$.
In \cite[Prop.
7.2]{LP2}, one finds the isomorphism
\begin{equation}\label{LPisomVstarfibdet} V_2^*(\O(\Delta)) \cong
(\O(p\Delta_1))_{|U}\end{equation} where $U$ is the complementary open subset
 of the base locus $\mathcal{B}$ in $S_{X_1}$.

\begin{proposition}\label{J[2]equivliftV} If $X$ is either a genus 2 curve or a non
hyperelliptic genus 3 curve, there is a $J[2]$-equivariant lifting
$\widetilde{V} : |2\Theta_1| \dashrightarrow |2\Theta|$ of the
rational map $V_2 : S_{X_1} \dashrightarrow S_X$.
\end{proposition}
\begin{proof} The
proposition is equivalent to the existence of a $J[2]$-equivariant
rational map
$\widetilde{V} : |2\Theta_1|^* \dashrightarrow |2\Theta|^*$
 such that the
following diagram commutes
\begin{equation*}\label{Vfactorsdiag}
\unitlength=0.6cm
\begin{picture}(10,4)
\put(1.3,3){$S_{X_1}$}  \put(8.3,3){$S_{X}$}
\put(1.2,0.5){$|2\Theta_1|^*$}  \put(8.2,0.5){$|2\Theta|^*$}
\put(7.6,3.2){\vector(1,0){0.1}} \dashline{0.2}(3,3.2)(7.6,3.2)
\put(7.6,0.7){\vector(1,0){0.1}} \dashline{0.2}(3,0.7)(7.6,0.7)
\put(1.9,2.8){\vector(0,-1){1.5}} \put(8.7,2.8){\vector(0,-1){1.5}}
\footnotesize{\put(0.9,1.8){$\varphi_{\Delta_1}$}
\put(8.9,1.8){$\varphi_{\Delta}$} \put(5.2,3.4){$V_2$}
\put(5.2,0.9){$\widetilde{V}$}}
\end{picture} \end{equation*}
In other words, one has to find a
factorization
\begin{equation}\label{Vfactors} W \xrightarrow
{\widetilde{V}^*} {\rm Sym}^p\, W_1
\xrightarrow{\varphi_{\Delta_1}^*}
H^0(S_{X_1},\,\O(p\,\Delta_1))\end{equation} of $V_2^* : W \to H^0(S_{X_1},\,\O(p\Delta_1))$, where the
first arrow defines a $J[2]$-equivariant map $\widetilde{V} :
|2\Theta_1| \dashrightarrow |2\Theta|^*$, and where the last
coincides with the canonical evaluation map.

Let us first prove the existence of a not necessarily $J[2]$-equivariant map $\widetilde{V}$.

\begin{lemma}\label{widetildeVexists} \textup{(1)} If $X$ has genus 2,
there is a rational map $\widetilde{V}$ as in diagram
\eqref{Vfactorsdiag}. It is defined by 4 homogeneous degree $p$
polynomials uniquely determined.\\
\textup{(2)} If $X$ is a general genus 3 curve (it is in particular
non hyperelliptic), the moduli space $S_X$ is projectively normal in
$|2\Theta|$. Therefore, there is a rational map $\widetilde{V}$ as
in diagram \eqref{Vfactorsdiag}. It is defined by 8 homogeneous
degree $p$ polynomials, uniquely determined modulo the Coble's
quartic.\end{lemma}

\begin{proof}(1) The base locus $\mathcal{B}$ is known to be finite (see
\cite{LnP}) and
\cite{Os1}) hence of
has codimension $\geq 2$ and the isomorphism
\eqref{LPisomVstarfibdet} extends to an isomorphism
$V_2^*(\O(\Delta)) \cong \O(p\Delta_1)$. Because the morphism
$\varphi_\Delta : S_X \to |2\Theta|^*$ is an isomorphism in that
case, there is a map $\widetilde{V}$ such that the diagram
\eqref{Vfactorsdiag} commute and it is uniquely determined.

(2) The work of Mochizuki (see \cite{Mo} and \cite[Thm 4.4]{Os2})
ensures that the base locus $\mathcal{B}$ has dimension 2. Because
$S_X$ has dimension 6 in that case, the isomorphism
\eqref{LPisomVstarfibdet}  extends again to an isomorphism
$V_2^*(\O(\Delta)) \cong \O(p\Delta_1)$. Because $S_{X}$ (resp. $S_{X_1}$ is normal and complete intersection in
$|2\Theta|^*$ (resp. $|2\Theta_1|^*$), it is projectively normal in $|2\Theta|^*$ (resp. $|2\Theta_1|^*$) (see
\cite[II, Exercise 8.4.(b), p. 188]{Ha}). Letting
$\O_{|2\Theta_1|^*}(n)$ denote the $n$-th power of the canonical
twisting sheaf on $|2\Theta_1|^*$, the canonical evaluation map,
which coincides with
\begin{equation*} H^0(|2\Theta_1|^*,\, \O_{|2\Theta_1|^*}(p))\cong {\rm Sym}^p
W_1 \xrightarrow{\varphi_{\Delta_1}^*}
H^0(S_{X_1},\,\O(p\Delta_1))\end{equation*} is therefore surjective
with kernel isomorphic to the image of
\begin{equation*}
H^0(|2\Theta_1|^*,\, \O_{|2\Theta_1|^*}(p-4))\cong {\rm Sym}^{p-4}
W_1 \xrightarrow {C_{X_1}} H^0(|2\Theta_1|^*,\,
\O_{|2\Theta_1|^*}(p))\cong {\rm Sym}^{p} W_1 \qedhere
\end{equation*}
\end{proof}

\begin{remark} \label{uniquenessg3p3} In the case of a (general) genus 3 curve, and in
characteristic 3 (the only case in which we will perform the
computations), the rational map $\widetilde{V}$ is uniquely
determined.\end{remark}

Because $J[2]$ does not act on sections but on multiplicative classes, we need to define actions of the corresponding Heisenberg group $\mathcal{G}(2)$.

\begin{lemma}\label{weightp^2action}
There is a Theta group homomorphism $\mathcal{G}(2) \to
\mathcal{G}_1(2)$ of weight $p$. In the cases of the previous lemma, it induces a weight $p^2$ action of
$\mathcal{G}(2)$ on both ${\rm Sym}^p\,W_1$ and
$H^0(S_{X_1},\,\O(p\,\Delta_1))$, compatible with the evaluation map.
\end{lemma}
\begin{proof} We check that the definition of \cite[p. 310]{Mu2} can be generalized in our case. Namely, there is a homomorphism of Heisenberg groups
$\eta_p
: {\mathcal{G}(\O(2p\,\Theta_1))}\to
{\mathcal{G}_1(2)}$ mapping any
$\gamma : \O(2p\,\Theta_1) \xrightarrow \sim T_x^* \O(2p\,\Theta_1)$
in $\mathcal{G}(\O(2p\Theta_1))$, $\eta_p(\gamma)$ to  the unique isomorphism
$\rho : \O(2\Theta_1) \xrightarrow{\sim}
T_{p\,x}^*\O(2\Theta_1)$ such that the composite
$$[p\,]^*\O(2\Theta_1) \xrightarrow \sim  \O(2p^2\,\Theta_1) \xrightarrow {\gamma^{\tens p}} T_x^*\O(2p^2\,\Theta_1) \xrightarrow \sim  [p\,]^*T_{p\,x}^*\O(2\Theta_1)$$ coincides with $[p\,]^*\rho$.
It raises elements of the center to the $p$-th power and the composite homomorphism $\eta_p \circ V^* : \mathcal{G}(2) \to
\mathcal{G}_1(2)$ has therefore weight $p$.

Using $p$-symmetric power, $\mathcal{G}_1(2)$ has a natural weight
$p$ action on ${\rm Sym}^p\,W_1$ and composing, it gives a weight $p^2$ action of $\mathcal{G}(2)$ on ${\rm Sym}^p\,W_1$. Because the evaluation map ${\rm Sym}^p\,W_1 \to H^0(S_{X_1},\,\O(p\,\Delta_1))$ is an isomorphism in genus 2 and has kernel invariant under the action of $\mathcal{G}(2)$ in case of a general genus 3 curve, it induces a weight $p^2$ action of $\mathcal{G}(2)$ on $H^0(S_{X_1},\,\O(p\,\Delta_1))$.\end{proof}

Because $\widetilde{V}^*$ is linear, there is no chance that it can be $\mathcal{G}(2)$-equivariant. We therefore define the subgroup $\mathcal{G}(2)[2]$ of order 2
elements in $\mathcal{G}(2)$ and
consider the induced actions on $W$, ${\rm Sym}^p\, W_1$ and
$H^0(S_{X_1},\,\O(p\Delta_1))$. Because $V_2^* : W \to
H^0(S_{X_1},\,\O(p\Delta_1))$ is linear and comes from the $J[2]$-equivariant map $V_2 :
S_{X_1} \dashrightarrow S_X$, it
is $\mathcal{G}(2)[2]$-equivariant.

Choose a theta structure on $
\mathcal{G}(2)$ and let $\{X_h,\,h\in H\}$ denote the associated theta
basis. Assume that the map $\widetilde{V}$ provided by the Lemma above is not $J[2]$-equivariant. Despite
$\widetilde{V}^*(X_0)$ might not be $\hat{H}$-invariant, its class
in $H^0(S_{X_1},\,\O(p\Delta_1))$ is for $V_2^*$ is $\mathcal{G}(2)[2]$-equivariant. The element
$$V_0':=\Frac{1}{2^g}\sum_{h^*\in \hat{H}}
h^*(\widetilde{V}^*(X_0))$$ in $({\rm Sym}^p\,W_1)^{\hat{H}}$ is
$\hat{H}$-invariant and maps onto $V_2^*(X_0)$. Defining
$V'_h:=h.V_0'$, the map $W \to {\rm Sym}^p W_1$ defined by $X_h \mapsto V_h'$ is $\mathcal{G}(2)[2]$-equivariant and induces $V_2^*$. The associated map $V' :
|2\Theta_1|^* \dashrightarrow |2\Theta|^*$ is
$J[2]$-equivariant  and we have proved the proposition.
\end{proof}

\begin{remark}\label{Ftheta} One checks that
$[p\,]^*i^*\gamma=V^*(F^*i^*)\gamma=V^*(\gamma^{\tens
p})=(V^*\,\gamma)^{\tens p}$ for any $\gamma$ in $\mathcal{G}(\O(2p\,\Theta_1))$. Therefore, the
homomorphism $\eta_p \circ V^*$  coincides with the homomorphism
$i^* : \mathcal{G}(2) \to \mathcal{G}_1(2)$ induced by the pull-back
by the quasi-isomorphism $i$. In particular, for any
two elements $\bar{\alpha}$ and $\bar{\beta}$ in $J[2]$,
\begin{equation*}e_{2,\,J_1}(F(\bar{\alpha}),\,F(\bar{\beta}))=e_{2,\,J}(\bar{\alpha},\,\bar{\beta})^p\end{equation*}
Because $J[2]$ is reduced and because $e_2$ takes its values in
$\mu_2$, we find that $F$ (hence $V$) is a symplectic isomorphism.
This implies that the choice of a Göpel system for $J[2]$ (resp. a theta structure on $\mathcal{G}(2)$) determines a Göpel system for $J_1[2]$ (resp. a theta structure $\mathcal{H}_1 \xrightarrow \sim \mathcal{G}_1(2)$
(where $\mathcal{H}_1: = \mathcal{H}\tens _{F_{\rm abs},\,k} k$)),
and that the associated  theta bases $\{X_\alpha\}_{\alpha\in \,H}$ and
$\{Y_{\alpha_1}\}_{\alpha_1 \in \,H_1}$ are compatible in the sense that
$Y_{\alpha_1}=i^*X_{V(\alpha_1)}$.
\end{remark}

\subsection{Frobenius action and Prym varieties}

The functoriality of Frobenius is also compatible with the correspondence between order 2 line bundles over $X$ and double étale covers of $X$.
If $F$ is the relative Frobenius, \cite[I.11]{SGA1} says that the
diagram
\begin{equation}\label{SGA}
\unitlength=0.6cm
\begin{picture}(8,3)
\put(1,2.5){$\widetilde{X}$} \put(6,2.5){$\widetilde{X}_1$}
\put(1,0){$X$}\put(6,0){$X_1$} \put(1.6,2.7){\vector(1,0){4}}
\put(1.6,0.2){\vector(1,0){4}} \put(1.2,2.3){\vector(0,-1){1.5}}
\put(6.2,2.3){\vector(0,-1){1.5}} \small \put(3.5,0.4){$F$}
\put(3.5,2.9){$F$} \put(0.5,1.4){$\pi$} \put(6.6,1.4){$\pi_1$}
\end{picture}\end{equation}
is cartesian.
As a consequence, the morphisms $\pi^* : J \to \widetilde{J} $ and ${\rm Nm} :
\widetilde{J} \to J$ commute with $V$.
\begin{proposition}\label{Prymordinary}
The following diagram
\begin{center}
\unitlength=0.6cm
\begin{picture}(10,3.5)
\put(0.3,0){\put(0.8,0){\put(0.2,2.5){$J\times P_\tau$}
 \put(0,0){$J_1\times
P_{\tau_1}$}\put(1.1,0.7){\vector(0,1){1.5}}

\put(0.8,0){\put(6,0){$\widetilde{J}_1$}
\put(6,2.5){$\widetilde{J}$} \put(1.6,2.7){\vector(1,0){4}}
\put(1.6,0.2){\vector(1,0){4}} \put(6.2,0.7){\vector(0,1){1.5}}}}

\small \put(4.7,2.9){$\sigma$} \put(4.7,0.4){$\sigma_1$}
\put(8,1.1){$V$}} \small \put(0,1.1){$V \times V_\tau$}
\end{picture} \end{center} is commutative.\\
Furthermore, $\sigma$ induces an isomorphism $J[p\,]\times
P_\tau[p\,] \xrightarrow{\sim } \widetilde{J}[p\,]$. In particular,
if $J$ is an ordinary abelian variety, then $\widetilde{J}$ is
ordinary if and only if $P_\tau$ is ordinary.
\end{proposition}
\begin{proof} The
Prym variety $$P_{\tau_1}:=\text{ker}(\text{Nm})^0 \subseteq
\widetilde{J}_1$$ coincides with the $p$-twist of
$P_\tau:=\text{ker}(\text{Nm})^0 \subseteq \widetilde{J}$ and it is
mapped by $V$ to $P_\tau$. The restriction
$V_{|P_{\tau_1}} : P_{\tau_1} \to P_\tau$ coincides with the Verschiebung $V_\tau : P_{\tau_1}
\to P_\tau$ for $P_\tau$ and the commutation of the diagram
follows from the
the fact that $V$ is a homomorphism and commutes with $\pi^*$.\\
The isomorphism $J[p\,]\times
P_\tau[p\,] \xrightarrow{\sim } \widetilde{J}[p\,]$ follows from the inclusion $\text{ker}\, \sigma \subseteq J[2] \times P_\tau[2]$
and the hypothesis $p\geq 3$. Because $X$ was supposed to be
ordinary, the last assertion comes from the fact that the ordinariness of any abelian variety can be read on the reduced part of its $p$-torsion. \end{proof}

For later use, let us mention the following result due to Nakajima,
proving that for a sufficiently general curve $X$, all the abelian
varieties appearing in the Proposition above are ordinary (one can
find a proof in \cite{Zh}).

\begin{proposition}\label{Nakajima} Let $X$ be a general, proper and smooth connected curve
over an algebraically closed field of characteristic $p$ and let $f
: Y \to X$ be an étale cover with abelian Galois group G. Then $Y$
is ordinary.\end{proposition}

Choose an element $z$ in $S_{\tau}=\{z\in J|\,z^2=\tau \} \subset
J[4]$ and let $z_1$ be $F(z)$ (equivalently $i^*z$). Let $\nu_p$ be
0 if $p\equiv 1$ mod 4 and be 1 if $p\equiv 3$ mod 4, then
\begin{equation*}
F^*z_1=V(z_1)=(-1)^{\frac{p-1}{2}}\, z=z\tens
\tau^{\nu_p}\end{equation*} Therefore, for any $L_1$ of
$P_{\tau_1}$, the cartesian square \eqref{SGA} gives
\begin{equation*}F^*((\pi_1)_*(L_1) \tens z_1)\cong \pi_*(F^*(L_1)) \tens
F^*z_1 \cong \pi_*(F^*(L_1)) \tens z \tens \tau^{\nu_p} \cong
\pi_*(F^*(L_1)) \tens z \qedhere
\end{equation*}
Assuming that the curve $X$ is
sufficiently general, the Proposition \ref{Nakajima} ensures that
all the Prym varieties are ordinary and one can choose symmetric principal divisors $\Xi$ and $\Xi_1$ on $P_\tau$ and $P_{\tau_1}$ respectively such that $\O(p\,\Xi_1) \cong V_\tau^*\O(\Xi)$ \cite[Lemma
2.2]{LP1}.

Let $\varphi_{2\Xi}$ be the
canonical map $P_\tau \twoheadrightarrow P_\tau/\pm \subseteq
|2\Xi|^*$ and define $\varphi_{2\Xi_1}$ analogously. Let $V_\tau^\pm : P_{\tau_1}/\pm \to P_{\tau}/\pm$ be the morphism induced by $V_\tau : P_{\tau_1} \to
P_\tau$ (which commutes with $[-1]$).
Assume that, as in the diagram \eqref{Vfactorsdiag}, there is a
$J[2]$-equivariant rational map $\widetilde{V} : |2\Theta_1|
\dashrightarrow |2\Theta|$ lifting $V_2 : S_{X_1} \dashrightarrow
S_X$. Restricting the map $\widetilde{V}$ to
$|2\Xi_1|^*$ (which identifies with one of the two connected components of $|2\Theta_1|^{\tau_1}$) yields a rational map $\widetilde{V}_{|2\Xi_1|^*} :
|2\Xi_1|^* \dashrightarrow |2\Xi|^*$ that lifts the morphism
$V_\tau^\pm : P_{\tau_1}/\pm \to P_{\tau}/\pm$. Therefore, one
obtains the following commutative diagram
\begin{equation} \label{restrictingVtildtoPrymdiag}
\unitlength=0.6cm
\begin{picture}(20,4)
\put(1,3){$P_\tau$}  \put(0.9,0){$P_{\tau_1}$}
\put(5.3,3){$P_{\tau}/\pm$} \put(5.2,0){$P_{\tau_1}/\pm$}
\put(10,3){$|2\Xi|^*$} \put(9.8,0){$|2\Xi_1|^*$}
\put(14,3){$|2\Theta|^\tau$} \put(13.7,0){$|2\Theta_1|^{\tau_1}$}
\put(18,3){$|2\Theta|$} \put(17.8,0){$|2\Theta_1|$}

\put(2,3.2){\vector(1,0){2.5}} \put(8,3){$\subseteq $}
\put(12.3,3){$\subseteq $}  \put(16.5,3){$\subseteq $}
\put(2,0.2){\vector(1,0){2.5}} \put(8,0){$\subseteq $}
\put(12.3,0){$\subseteq $}  \put(16.5,0){$\subseteq $}

\put(1.3,1){\vector(0,1){1.6}} \put(6,1){\vector(0,1){1.6}}
\dashline{0.2}(10.6,1)(10.6,2.4)\put(10.6,2.4){\vector(0,1){0.2}}
\dashline{0.2}(18.6,1)(18.6,2.4)\put(18.6,2.4){\vector(0,1){0.2}}

\small \put(3,3.5){$\varphi_{2\Xi}$} \put(3,0.5){$\varphi_{2\Xi_1}$}
\put(0.4,1.5){$V_\tau$} \put(5,1.5){$V_\tau^\pm$}
\put(11,1.5){$\widetilde{V}_{|2\Theta_{1}|^*}$}
\put(19,1.5){$\widetilde{V}$}
\end{picture} \end{equation}
In terms of coordinates functions, we have therefore proven the
following :

\begin{proposition} \label{widetildeVJ[2]equivariant} Assume that $X$ is proper
and smooth genus $g=2$ or 3 curve over $k$, sufficiently general.
Then, for any non zero $\tau$ in $J[2]$,  there is a
$P_\tau[2]$-equivariant map $\widetilde{V}_\tau : |2\Xi_1|^*
\dashrightarrow |2\Xi|^*$ such that the pull-back $V_\tau^* :
H^0(P_\tau,\,\O(2\Xi)) \to H^0(P_{\tau_1},\,\O(2p\,\Xi_1))$ factors
as the composite
\begin{eqnarray}H^0(P_\tau,\,\O(2\Xi)) \xrightarrow{\widetilde{V}_\tau^*} {\rm
Sym}^p\,H^0(P_{\tau_1},\,\O(2\Xi_1)) \to
H^0(P_{\tau_1},\,\O(2p\,\Xi_1))\end{eqnarray} where the last arrow
is the
evaluation map.\\
 Letting $\wt$ (resp. $\wt_1$) be the order 2 lift
of $\tau$ in $\mathcal{G}(2)$ (resp. $\tau_1$ in $\mathcal{G}_1(2)$)
associated to $z$ (resp. $z_1$),
$\widetilde{V}_\tau$ agrees (up to a multiplicative scalar) with the
restriction
\begin{eqnarray*}\widetilde{V}_{|(W^{\wt})^*} :  \left(W^{\wt}\right)^* \longrightarrow {\rm
Sym}^p\,\left(W_1^{\wt_1}\right)^* \end{eqnarray*}
\end{proposition}

\section{The proof of the Theorem \ref{MainTheorem}} \label{proof}

From now on, we choose once for all a theta structure $\widetilde{\phi}_0
: \mathcal{H} \xrightarrow \sim \mathcal{G}(2)$ and we let $\{X_\alpha,\,\alpha \in H\}$ be the associated theta basis. It gives a system of homogeneous
coordinates $\{x_\alpha,\,\alpha \in H\}$ for $|2\Theta|$.
Recall from Remark \ref{Ftheta} that this choice also determines a theta structure
$\widetilde{\phi}_1 : \mathcal{H}_1 \xrightarrow \sim
\mathcal{G}_1(2)$, hence a theta basis for $W_1$ and coordinate
functions for $|2\Theta_1|$ adapted to the action of $J[2]$ that we denote $\{Y_\alpha|\,\alpha \in H\}$ and $\{y_\alpha|\,\alpha \in H\}$ respectively. Also, given
$\tau=(\alpha_0,\,\alpha_0^*)$ in $H\times \hat{H}$, recall that the choice of
a lift $\wt=(\mu,\alpha_0,\,\alpha_0^*)$ (with
$\mu^2=\alpha_0^*(\alpha_0)$) determines a lift
$\wt_1=(\mu^p,\alpha_0,\,\alpha_0^*)$ of $\tau_1$ in $\mathcal{H}_1$. Fiw once for all a square root $\mu_0$ of $-1$ in $k$. We will take $\mu=1$ if $\alpha_0^*(\alpha_0)=1$
and $\mu=\mu_0$ otherwise.

\subsection{From geometry to linear algebra}
The direct sum ${\rm Sym}^p\,(W_1^{\wt_1})^*\oplus {\rm
Sym}^p\,(W_1^{-\wt_1})^*$ (which depends only on the choice of
$\tau$) is endowed with an action of $\mathcal{H}$ of weight $p^2$
 and the quotient map
\begin{equation*}{\rm Sym}^p\,W_1^*\to {\rm Sym}^p\,(W_1^{\wt_1})^*\oplus
{\rm Sym}^p\,(W_1^{-\wt_1})^*\end{equation*} is
equivariant for the action of $\mathcal{H}$ on both spaces. Taking
all order 2 elements of $J[2]$ together, we find a morphism of
$\mathcal{H}$-representations \begin{equation*} R : {\rm
Sym}\,^p\,W_1^* \to S_p:=\bigoplus_{\tau \in\, J[2]\setminus
\{0\}}{\rm Sym}^p\,(W_1^{\wt_1})^*\oplus {\rm
Sym}^p\,(W_1^{-\wt_1})^*.\end{equation*} Because $\widetilde{V}$ is
given by a linear sub-space of ${\rm Sym}\,^p\,W_1^*$, isomorphic to
$W^*$ and endowed with an (irreducible) action of $\mathcal{G}(2)$
of weight $p^2$, it is determined by its $\bar{H}$-invariant part
and Theorem \ref{MainTheorem} follows from the Proposition

\begin{proposition} \label{hatRinjective} When $g=2$ or 3 and $p=3,\,5$ or 7, the restriction map
\begin{equation*} \hat{R} : \left({\rm
Sym}\,^p\,W_1^*\right)^{\hat{H}} \to S_p\end{equation*} is
injective.\end{proposition}

\begin{remark}\label{methodlimited} Certainly, a necessary condition
is that \begin{equation*} \dim {\rm
Sym}^p\,W_1^*=\left(\begin{array}{c}2^g+p-1\\
2^g-1\end{array}\right)\leq \dim S_p=2(2^{2g}-1)\left(\begin{array}{c}2^{g-1}+p-1\\
2^{g-1}-1\end{array}\right).\end{equation*} This cannot be the case for
large $p$. More precisely, when $g=2$, $\dim\, {\rm Sym}^p\, W_1^*>\dim \,S_p$
for $p>7$, and when $g=3$, $\dim\, {\rm Sym}^p\, W_1^*>\dim \,S_p$
for $p>11$.\end{remark}

\subsection{Preparation in small genus}

 Since ${\rm Sym}\,^p\,W_1^*$ is
generated by the free family of monomials $\prod_{\alpha\in \,H}
y_{\alpha}^{e_{\alpha}}$ with $\sum e_\alpha=p$, $\left({\rm
Sym}\,^p\,W_1^*\right)^{\hat{H}}$ is generated by the subfamily
consisting  of elements whose set of exponents satisfies
\begin{equation*}\sum_{\alpha \in H|\,\alpha_0^*(\alpha)=-1}e_\alpha
\equiv 0 \ \text{mod. }2 \text{ for all }\alpha_0^*\text{ in
}\hat{H}.\end{equation*} Writing such an element under the form
$\prod_{\alpha\in \,H}
y_{\alpha}^{\bar{e}_{\alpha}}\left[\prod_{\alpha\in \,H}
y_{\alpha}^{f_{\alpha}}\right]^2$ with $\bar{e}_\alpha=0$ or 1 and
$\sum \bar{e}_\alpha+2\sum f_\alpha=p$, we find that
\begin{equation*}\sum_{\alpha \in H|\,\bar{e}_\alpha=1}\,\alpha=0.\end{equation*}
It is easily seen that $E$ can be $\{0\}$ and $H-\{0\}$. In genus 2,
these are the only possibilities. In genus 3, one also has  to
consider the case where $E$ has cardinal 3 and $E'=\{0\}\cup E$ is a
cardinal 4 subgroup of $H$ (the set of such subgroups is in 1-1
correspondance with the set of non zero elements of $\hat{H}$).
Consider therefore the sets
\begin{equation*}\label{Afamily}A(p)=\left\{A_{\underline{f}}=y_{0}\left[\prod_{\alpha
\in \,H} y_{\alpha}^{f_{\alpha}}\right]^2,\text{with
}|\underline{f}|=\Frac{p-1}{2}\right\}\end{equation*}
\begin{equation*}\label{Bfamily}B(p)=\left\{B_{\underline{f}}=\prod_{\alpha \in H-\{0\}}
y_\alpha \left[\prod_{\alpha\in \,H}
y_{\alpha}^{f_{\alpha}}\right]^2,\text{with
}|\underline{f}|=\Frac{p-2^g+1}{2}\right\}\end{equation*}and
\begin{equation*}\label{Cfamily}C(p)=\left\{C_{\alpha^*,\,\underline{f}}=\prod_{\alpha \in H-\{0\},\,\alpha^*(\alpha)=1}
y_\alpha \left[\prod_{\alpha\in \,H}
y_{\alpha}^{f_{\alpha}}\right]^2,\text{with }\alpha^* \text{ in
}\hat{H}-\{0\} \text{ and }
|\underline{f}|=\Frac{p-3}{2}\right\}\end{equation*}
 (where $\underline{f}$ is, in each case, the multi-index
 $(f_\alpha,\,\alpha \in H)$ with $|\underline{f}|= \sum
f_\alpha$).

To compute the image of these monomials in the various ${\rm
Sym}^p\,(W_1^{\wt_1})^*\oplus{\rm Sym}^p\,(W_1^{-\wt_1})^*$, it is
convenient to distinguish whether $\tau$ is an element of $\hat{H}$
or not. In the first case, the following lemma
is straightforward
\begin{lemma} \label{reductionalpha0=0} For a given
$\tau=\wt =\alpha_0^*$, all the monomials in $A(p),\ B(p)$ or $C(p)$
map to 0 in ${\rm Sym}^p\,(W_1^{-\wt_1})^*$. The only monomials in
$A(p)\bigsqcup B(p)$ (resp. $A(p)\bigsqcup B(p) \bigsqcup C(p)$) not
mapping to 0 in ${\rm Sym}^p\,(W_1^{\wt_1})^*$ are those that can be
written under the form
\begin{equation*} A_{\underline{f}}=y_{0}\left[\prod_{\alpha \in \,H,\, \alpha_0^*(\alpha)=1}
y_{\alpha}^{f_{\alpha}}\right]^2
\end{equation*}
(resp. as well as
\begin{equation*} \left.C_{\alpha_0^*,\,\underline{f}}=\prod_{\alpha \in H-\{0\},\,\alpha_0^*(\alpha)=1}
y_\alpha \left[\prod_{\alpha \in \,H}
y_{\alpha}^{f_{\alpha}}\right]^2\right)
\end{equation*}
and they all map to a different monomial in ${\rm Sym}^p
\,(W_1^{\wt_1})^*$.\end{lemma}

\subsection{Proof of Proposition \ref{hatRinjective} in the genus 2 case}
As indicated above, a basis of $\left({\rm
Sym}\,^p\,W_1^*\right)^{\hat{H}}$ in genus 2 is the (disjoint) union
of the two sets $A(p)$ and $B(p)$.

\begin{lemma} \label{reductionalpha0neq0g2} Assume that
\begin{equation}\label{V0shapeg2} V_0=\sum_{A(p)} a_{\underline{f}} A_{\underline{f}}+
\sum_{B(p)} b_{\underline{f}} B_{\underline{f}}\end{equation} maps
to zero in $S_p$. Then, for any integer $0\leq k\leq \frac{p-1}{2}$
\begin{equation*}\sum_{f_{00}+f_{01}=k\text{ and }f_{01}+f_{11} \text{ even
}}a_{\underline{f}}=0;\hspace{1cm} \sum_{f_{00}+f_{01}=k\text{ and
}f_{01}+f_{11} \text{ odd }}a_{\underline{f}}=0\end{equation*}
\begin{equation*}
\sum_{f_{00}+f_{10}=k\text{ and }f_{10}+f_{11} \text{ even
}}a_{\underline{f}}=0;\hspace{1cm} \sum_{f_{00}+f_{10}=k\text{ and
}f_{10}+f_{11} \text{ odd }}a_{\underline{f}}=0\end{equation*}
\begin{equation*}\sum_{f_{00}+f_{11}=k\text{ and }f_{10}+f_{11} \text{ even
}}a_{\underline{f}}=0;\hspace{1cm} \sum_{f_{00}+f_{11}=k\text{ and
}f_{10}+f_{11} \text{ odd }}a_{\underline{f}}=0\end{equation*} and
we have analogous equalities for the $b_{\underline{f}}$ (with
$0\leq k \leq \frac{p-3}{2}$).
\end{lemma}
\begin{proof} Choose a non zero $\alpha_0$ in $H$ and for
$\tau=(\alpha_0,\,\alpha_0^*)$, let $\wt$ be defined as above. In
particular, $\wt_1=(\mu^p,\,\alpha_0,\,\alpha_0^*)$. A generating
system for $W_1^{\wt_1}$ (resp. $W_1^{-\wt_1}$) is
$$\{(Y_\alpha+\wt_1.Y_\alpha),\,(\alpha \in H)\} \text{ (resp.
}\{(Y_\alpha-\wt_1.Y_\alpha),\,(\alpha \in H)\})$$Since one has
$\wt_1.Y_\alpha=\mu^p \alpha_0^*(\alpha)Y_{\alpha+\alpha_0}$,
 one finds, letting $\bar{y}_\alpha$ denote the class of $y_\alpha$
 in $(W_1^{\wt_1})^*$, that $\bar{y}_{\alpha+\alpha_0}=\mu^{p}\alpha_0^*(\alpha)
\bar{y}_{\alpha}$. Similarly, letting $\bar{y}_\alpha$ denote also
the class of $y_\alpha$
 in $(W_1^{-\wt_1})^*$ (there will be no risk of confusion), one finds that $\bar{y}_\alpha=-\mu^{p}\alpha_0^*(\alpha)
\bar{y}_{\alpha}$.

Choose a non zero $\alpha$ in $H$ such that $H=<\alpha_0,\,\alpha>$.
Then $A_{\underline{f}}$ maps to
\begin{equation}\label{AftaunotinhatHg2}\alpha_0^*(\alpha_0)^{f_{\alpha_0}+f_{\alpha+\alpha_0}} \left(
\bar{y}_0^{1+2(f_0+f_{\alpha_0})}
{\bar{y}_\alpha}^{2(f_\alpha+f_{\alpha+\alpha_0})}\right)\end{equation}in
both  ${\rm Sym}^p\,(W_1^{\wt_1})^*$ and ${\rm
Sym}^p\,(W_1^{-\wt_1})^*$. Similarly, the monomial
$B_{\underline{f}}$ maps to
\begin{equation}\label{BftaunotinhatHg2}
\alpha_0^*(\alpha).\alpha_0^*(\alpha_0)^{f_{\alpha_0}+f_{\alpha+\alpha_0}}
\left( \bar{y}_0^{1+2(f_0+f_{\alpha_0})}
{\bar{y}_\alpha}^{2(f_\alpha+f_{\alpha+\alpha_0})}\right)\end{equation}
in both  ${\rm Sym}^p\,(W_1^{\wt_1})^*$ and ${\rm
Sym}^p\,(W_1^{-\wt_1})^*$.

Specifically, fix $\alpha_0=01$, choose $\alpha=10$ and a positive integer
$k$, and look at the coefficient of the monomial
$\bar{y}_{00}^{1+2k}\bar{y}_{10}^{p-1-2k}$ in the image of $V_0$ in
${\rm Sym}^p\,(W_1^{\wt_1})^*$ for all the elements
$\tau=(01,\,\alpha_0^*)$ in $J[2]$. Using \eqref{AftaunotinhatHg2},
\eqref{BftaunotinhatHg2} above, we find the following expressions
\begin{equation*}\begin{array}{c l}
(\alpha_0^*=00) : & \sum_{f_{00}+f_{01}=k} a_{\underline{f}}+
\sum_{f_{00}+f_{01}=k} b_{\underline{f}}\vspace{0.3cm}\\
(\alpha_0^*=01) : & \sum_{f_{00}+f_{01}=k}
(-1)^{f_{01}+f_{11}}a_{\underline{f}}-
\sum_{f_{00}+f_{01}=k} (-1)^{f_{01}+f_{11}} b_{\underline{f}}\vspace{0.3cm} \\
(\alpha_0^*=10) : & \sum_{f_{00}+f_{01}=k} a_{\underline{f}}-
\sum_{f_{00}+f_{01}=k} b_{\underline{f}}\vspace{0.3cm} \\
(\alpha_0^*=11) : & \sum_{f_{00}+f_{01}=k}
(-1)^{f_{01}+f_{11}}a_{\underline{f}}+ \sum_{f_{00}+f_{01}=k}
(-1)^{f_{01}+f_{11}} b_{\underline{f}}\end{array}\end{equation*} If
$\hat{R}(V_0)=0$, we easily derive from the expressions above the first line in the equations of the statement for both the $a_{\underline{f}}$'s and the $b_{\underline{f}}$'s. Repeating the process with $\alpha_0=10$ or $11$, $\alpha=01$ in
both cases, and letting $k$ vary through the relevant sets proves
the lemma. \end{proof}

The proof of proposition \ref{hatRinjective} in genus 2 now boils
down, using Lemmas \ref{reductionalpha0=0} and
\ref{reductionalpha0neq0g2}, to an easy exercise of linear algebra
left to the reader.

\subsection{Proof of Proposition \ref{hatRinjective} in the genus 3 case}
In genus 3, a basis of $\left({\rm Sym}\,^p\,W_1^*\right)^{\hat{H}}$
in the (disjoint) union of the three sets $A(p)$, $B(p)$ and $C(p)$.
The following lemma is very similar to the
analogous Lemma in genus 2 case, though more complicated and more
inconvenient to write down. We leave the proof to the reader and we simply indicate the tricks we found useful to do the computations. For any non zero $\alpha_0$ in $H$, we will choose a
subgroup $H(\alpha_0)$ in $H$, not containing $\alpha_0$, such that
$H(\alpha_0)$ and $\alpha_0$ together generate $H$ (it is analogous
to the choice of an element $\alpha$ in the proof of Lemma
\ref{reductionalpha0neq0g2}). Notice that $H(\alpha_0)$ inherits an
ordering from the lexicographic order on $H$. More specifically, we
will choose $H(001)=\{000,\,010,\,100,\,110\}$,
$H(010)=H(011)=\{000,\,001,\,101,\,101\}$, and
$H(\alpha_0)=\{000,\,001,\,010,,011\}$ otherwise. For any
multi-index $\underline{f}$, let $[\underline{f}]_{\alpha_0}$ denote
the 4-tuple $(f_\alpha+f_{\alpha_0+\alpha})_{\alpha \in
H(\alpha_0)}$ (with lexicographic order). Also, let
$\Sigma_{\alpha_0}\underline{f}$ stand for the sum $\sum_{\alpha \in
H(\alpha_0)} f_{\alpha_0+\alpha}$. We will keep on denoting by
$|\underline{f}|$ the sum $\sum f_{\alpha}$ for any $n$-tuple.

\begin{lemma} \label{reductionalpha0neq0g3} Assume that
\begin{equation*}V_0=\sum_{A(p)} a_{\underline{f}} A_{\underline{f}}+
\sum_{B(p)} b_{\underline{f}}
B_{\underline{f}}+\sum_{C(p)}c_{\alpha^*,\,\underline{f}}
C_{\alpha^*,\,\underline{f}}\end{equation*} maps to zero in $S_p$.
Then, for all 4-tuples
$\underline{k}=(k_{00},\,k_{01},\,k_{10},\,k_{11})$ of positive
integers
 such that
$|\underline{k}|=\frac{p-1}{2}$ and for all non-zero $\alpha_0$ in $H$, we have the equalities
\begin{equation*} \sum
_{[\underline{f}]_{\alpha_0}=\underline{k}\text{ and }
\Sigma_{\alpha_0}(\underline{f}) \text{ even}} a_{\underline{f}}+
\sum
_{[\underline{f}]_{\alpha_0}=(k_{00},\,k_{01}-1,\,k_{10}-1,\,k_{11}-1)
\text{ and } \Sigma_{\alpha_0}(\underline{f}) \text{ even}}
b_{\underline{f}}=0
\end{equation*}
and \begin{equation*} \sum
_{[\underline{f}]_{\alpha_0}=\underline{k}\text{ and }
\Sigma_{\alpha_0}(\underline{f}) \text{ odd}} a_{\underline{f}}+
\sum
_{[\underline{f}]_{\alpha_0}=(k_{00},\,k_{01}-1,\,k_{10}-1,\,k_{11}-1)
\text{ and } \Sigma_{\alpha_0}(\underline{f}) \text{ odd}}
b_{\underline{f}}=0
\end{equation*}
Also, for all 4-tuples
$\underline{k}=(k_{00},\,k_{01},\,k_{10},\,k_{11})$ of positive
integers
 such that
$|\underline{k}|=\frac{p-3}{2}$ and for all non-zero $\alpha_0$ in $H$, we have the equalities
\begin{equation*} \sum _{[\underline{f}]_{\alpha_0}=\underline{k}\text{ and }
\Sigma_{\alpha_0}(\underline{f}) \text{ even}}
c_{\alpha^*,\,\underline{f}}=0 \text{ and } \sum
_{[\underline{f}]_{\alpha_0}=\underline{k}\text{ and }
\Sigma_{\alpha_0}(\underline{f}) \text{ odd}}
c_{\alpha^*,\,\underline{f}}=0 \end{equation*} for all $\alpha^*$ in
$\hat{H}-\{0\}$.
\end{lemma}

Combining the data given by the Lemmas \ref{reductionalpha0=0} and
\ref{reductionalpha0neq0g3}, we reduce once again our problem to an
easy question of linear algebra. Nonetheless, it involves a large
number of unknowns and equations and those equations should not be
written bluntly since a few (easy) tricks simplify the matter
considerably.

\section{Elliptic curves, Kummer's quartic surface and Coble's quartic
hypersurface}\label{KumCob}

In this section, we give some preparation for the computations to
come in the next section.

\subsection{Kummer's quartic surface and associated elliptic curves}

Let us begin with recalling some well-known results dealing with the
geometry of the Kummer's quartic surface. Proofs can be found in
\cite{GD} (PhD thesis) or in  \cite{GH} in the complex case but they
can be carried over to any algebraically closed base field $k$ of
characteristic different from 2.

\begin{lemma} \label{Kummerdetails} {\rm (1)} Let $X$ be a smooth and projective genus 2 curve over $k$ and let $J$ be its Jacobian. The
scheme-theoretic image of the morphism ${\rm K}_X : J \to |2\Theta|$
identifies with the quotient of $J$ under the action of $\{\pm \}$.
It is a reduced, irreducible, $J[2]$-invariant quartic in
$|2\Theta|$ with 16 nodes and no other singularities, i.e., a
\emph{Kummer surface}.\\
\textup{(2)} In the coordinate system $\{x_\bullet\}$ defined above,
there are scalars $k_{00},\,k_{01},\,k_{10}$ and $k_{11}$ such that
the equation defining the Kummer quartic surface ${\rm Kum}_X$ is
\begin{equation}\label{eqKummer}K_X=S^K+ 2k_{00}
P^K+k_{01}Q^K_{01}+k_{10}Q^K_{10}+k_{11}Q^K_{11}\end{equation} where
$$\begin{array}{c}S^K= x_{00}^4+x_{01}^4+x_{10}^4+x_{11}^4,\hspace{1cm}
P^K=x_{00}x_{01}x_{10}x_{11},\\
Q^K_{01}=x_{00}^2x_{01}^2+x_{10}^2x_{11}^2,\hspace{1cm}Q^K_{10}=x_{00}^2x_{10}^2+x_{01}^2x_{11}^2,\hspace{1cm}
Q^K_{11}=x_{00}^2x_{11}^2+x_{01}^2x_{10}^2.\end{array}$$
\textup{(3)} These scalars $k_{00},\,k_{01},\,k_{10}$ and $k_{11}$
satisfy the cubic relationship
\begin{equation}\label{relKummercoef}4+k_{01}k_{10}k_{11}-k_{01}^2-k_{10}^2-k_{11}^2+k_{00}^2=0\end{equation}
and one has \begin{equation} \label{coefconstraints}
\left\{\begin{array}{l} k_{01}\neq \pm 2, \ k_{10}\neq \pm
2,\  k_{11}\neq \pm 2,\\
k_{01}+k_{10}+k_{11}+2\pm k_{00}\neq 0,\\ k_{01}+k_{10}-k_{11}-2\pm
k_{00} \neq 0,\\ k_{01}-k_{10}+k_{11}-2\pm k_{00} \neq 0,\\
-k_{01}+k_{10}+k_{11}-2\pm k_{00} \neq 0 \end{array}\right.
\end{equation}\end{lemma}

Let $\tau=(\alpha_0,\,\alpha_0^*)$ be a non zero element of
$J[2]=H\times \hat{H}$. Fix an order 2 lift $\wt$ of $\tau$ in $\mathcal{H}$ as in the previous section.
 The space $W$ splits in the direct sum
$W =W^{\wt}\oplus W^{-\wt}$ of the two 2-dimensional spaces of
eigenvectors
of $\wt$. Denote by $\Delta^+(\wt)$ (resp. $\Delta^-(\wt)$) the
corresponding projective lines in $|2\Theta|$. Again, one can find a non zero $\alpha=\alpha(\tau)$ in
$H$ such that the images $\bar{x}_0$ and
$\bar{x}_\alpha$ (of $x_0$ and $x_\alpha$ respectively) via the
canonical  map $W^* \twoheadrightarrow (W^{\wt})^*$ give a set of
homogeneous coordinates for $\Delta^+(\wt)$. We will write
$\lambda_0$ for $\bar{x}_0$ and $\lambda_1$ for $\bar{x}_\alpha$.

\begin{remark} \label{fancybases} We can similarly construct a system of coordinates
$\{\bar{\lambda}_0,\,\bar{\lambda}_1\}$ for $\Delta^-(\wt)$. Via Wirtinger's isomorphism, it gives dually a
basis
$\{\Lambda_0,\,\Lambda_1,\,\bar{\Lambda}_0,\,\bar{\Lambda}_1\}$ of
$W$ that splits into bases $\{\Lambda_0,\,\Lambda_1\}$ and
$\{\bar{\Lambda}_0,\,\bar{\Lambda}_1\}$ of $W^{\wt}$ and $W^{-\wt}$
respectively. One then notices that this is the theta basis associated to a suitable theta structure on $\mathcal{G}(2)$.
\end{remark}

Restricted to $\Delta^+(\wt)$, the equation of the Kummer surface
reduces, up to a multiplicative scalar,  to
\begin{equation}\label{KumcapDelta}
\lambda_0^4+\lambda_1^4+\omega \lambda_0^2\lambda_{1}^2=0
\end{equation} where $\omega:=\omega(\tau)$ depends on $\tau$ but not on
our particular choice of a lifting $\wt$ of $\tau$. In particular,
the equation
of the Kummer surface restricts in the same way to $\Delta^-(\wt)$.
The points
of $\Delta^+(\wt)\cap {\rm Kum}_X$ correspond to the four points of
$P_\tau[2]$ (in particular, $\omega\neq \pm 2$) and they must have homogeneous coordinates $(a :
\,b),\,(a : \,-b),\,(b : \,a)$ and $(b : \,-a)$, whence $\omega =-\Frac{a^4+b^4}{a^2b^2}$.
It is easy to compute the various $\omega(\tau)$ in terms of the
$k_\bullet$. Namely,

\begin{equation}\label{formomegatau} \begin{array}{|c|c|c|c|c|}

\hline \alpha_0 & \alpha_0^*=00 & \alpha_0^*=01 & \alpha_0^*=10 & \alpha_0^*=11 \\
\hline 00 & \star & k_{10}  & k_{01}& k_{11} \\
\hline 01 & \Frac{2(k_{00}+k_{10}+k_{11})}{2+k_{01}} &
\Frac{2(-k_{00}+k_{10}-k_{11})}{2-k_{01}}&
\Frac{2(-k_{00}+k_{10}+k_{11})}{2+k_{01}} &
\Frac{2(k_{00}+k_{10}-k_{11})}{2-k_{01}}\\
\hline 10 & \Frac{2(k_{00}+k_{01}+k_{11})}{2+k_{10}} &
\Frac{2(-k_{00}+k_{01}+k_{11})}{2+k_{10}}&
\Frac{2(-k_{00}+k_{01}-k_{11})}{2-k_{10}} &
\Frac{2(k_{00}+k_{01}-k_{11})}{2-k_{10}}\\
\hline 11 & \Frac{2(k_{00}+k_{01}+k_{10})}{2+k_{11}} &
\Frac{2(k_{00}+k_{01}-k_{10})}{2-k_{11}}&
\Frac{2(-k_{00}+k_{01}-k_{10})}{2-k_{11}} &
\Frac{2(-k_{00}+k_{01}+k_{10})}{2+k_{11}}\\
\hline
\end{array}\hspace{0.3cm} \end{equation} \\

Notice that the inequations \eqref{coefconstraints} ensure that
those coefficients are well-defined scalars, and they give another
reason why $\omega(\tau)$ cannot equal $\pm 2$. Because an elliptic curve $E$ is completely determined by the branch
locus of its Kummer map $E \to E/\{\pm\} \cong \P^1$, these data
allow one to determine the elliptic curve $P_\tau$ arising as the
Prym variety associated to the double cover corresponding to $\tau$.

\subsection{Coble's quartic and associated Kummer
surfaces}\label{sectCoblequartic}

Now, $X$ is a non hyperelliptic curve of genus 3 over $k$. Take
$\tau=(\alpha_0,\,\alpha_0^*)$ to be non-zero and fix a lift $\wt$
in $\mathcal{H}$.
 The space $W$ splits again in the direct sum
$W =W^{\wt}\oplus W^{-\wt}$ of the two 4-dimensional spaces of
eigenvectors $\wt$
and we let $\Delta^+(\wt)$ and  $\Delta^-(\wt)$ denote the corresponding
projective space $\P^3$ in $|2\Theta|$.
Again, one can find a cardinal 4 subgroup $H(\tau)$ of $H$ such that
the images $\bar{x}_\alpha$ of $x_\alpha$ ($\alpha$ in $H(\tau)$)
via the canonical map $W^* \twoheadrightarrow (W^{\wt})^*$ give a
set of homogeneous coordinates for $\Delta^+(\wt)$. This set of coordinates will be
denoted $\lambda_{00},\,\lambda_{01},\,\lambda_{10},\,\lambda_{11}$
in such a way that the lexicographical order is respected and if
$\alpha$ is in $H(\tau)$, we will write
$\bar{x}_\alpha=\lambda_{\bar{\alpha}}$ (with our conventions,
$\alpha \mapsto \bar{\alpha}$ is an group isomorphism). In
particular, $\bar{x}_0=\lambda_{\bar{0}}=\lambda_{00}$.

The following results are also well-known (the reader may refer to
\cite{Co}, \cite{GH} and \cite{Pa}).

\begin{lemma} There is a unique $J[2]$-invariant quartic ${\rm Cob}_X$ in
$|2\Theta|$ whose singular locus is the Kummer variety ${\rm Kum}_X$
associated to $X$.\\ If $\tau$ is a non zero element of $J[2]$ and
if $\Delta\cong \P^3$ is one of the two connected components of
$|2\Theta|^\tau$, then the intersection $\Delta\cap {\rm Cob}_X$ is
isomorphic to the Kummer surface ${\rm Kum}_\tau$ associated to the
genus 2 curve $Y_\tau$ whose Jacobian variety is isomorphic to  the
Prym variety $P_\tau$ corresponding to the double \'etale covering
$\widetilde{X} \to X$ associated to $\tau$.\\
The hypersurface ${\rm Cob}_X$ is completely determined by the data
of the 63 Kummer surface ${\rm Kum}_\tau$ defined above.
\end{lemma}

For later use, we give some details on the last point, involving
some calculations. Using the theta structure we have chosen, the equation of ${\rm Cob}_X$ goes under the
form
\begin{equation} \label{Coblequartic}
C_X =  S^C +\sum _{\alpha \in H-\{0\}} \gamma_\alpha Q^C_\alpha +
\sum _{\alpha^* \in \hat{H}
-\{0\}} \delta_{\alpha^*}
P^C_{\alpha^*}\end{equation} where, letting $\hat{\alpha}$ denote
the dual of $\alpha$ in $\hat{H}$,
$$S^C = \sum _{\beta \in H} x_{\beta}^4 ; \ \ Q^C_\alpha  =  \sum_{\beta \in H|\,
\hat{\alpha}(\beta)=1}x_\beta^2x_{\beta+\alpha}^2;  \ \
P^C_{\alpha^*}  =  \prod_{\beta \in H|\,
\alpha^*(\beta)=1}x_\beta+ \prod_{\beta \in H|\,
\alpha^*(\beta)=-1}x_{\beta}.
$$
Once again, the application $\alpha^*
\mapsto \{\beta \in H|,\,\alpha^*(\beta)=1\}$ gives a one-to-one correspondence between the set $\hat{H}-\{0\}$ and the set $H_4$
  of cardinal 4 subgroups of $H$. If $G$ is element of $H_4$, we denote by $\alpha_G^*$ the
corresponding element of $\hat{H}-\{0\}$ (in such a way that $G=\{\beta \in
H|,\,\alpha_G^*(\beta)=1\}$) and we define
$\delta_GP_G^C:=\delta_{\alpha^*_G}P^C_{\alpha^*_G}$.

Restricted to $\Delta^+(\wt)$, the equation $C_X$ of the Coble's
quartic hypersurface reduces, up to a multiplicative scalar, to the
equation of a Kummer surface whose coefficients depend only on
$\tau$ (and not on $\wt$). Those coefficients
$k_{00}(\tau),\,k_{01}(\tau),\,k_{10}(\tau),\,k_{11}(\tau)$ are
determined in terms of the $\gamma_\alpha$'s and the
$\delta_\alpha$'s. Namely, one finds that either
\begin{equation}\label{CobletoKumalpha0=0} 2k_{00}(\alpha_0^*)=\delta_{\alpha_0^*}\text{ and }
k_{\bar{\alpha}}(\alpha_0^*)=\gamma_\alpha \text{ for all non zero
}\alpha\text{ in }H(\tau)
\end{equation}
if $\alpha_0=0$ or
\begin{equation}\label{CobletoKumalpha0neq0} \begin{array}{c}  k_{00}(\tau)=2\Frac{\delta_{H(\tau)} +\sum_{G \in H_4 |\, G\cap H(\tau)
=<\alpha_1>,\, \alpha_G^*(\alpha_0)=-1} \alpha_0^*(\alpha+\alpha_0)
\delta_{G} }{2+\alpha^*_0(\alpha_0)\gamma_{\alpha_0}} \vspace{0.5cm}\\
\text{and }\  k_{\bar{\alpha}}(\tau)=\Frac{2
(\gamma_\alpha+\alpha^*_0(\alpha_0)
\gamma_{\alpha+\alpha_0})+\delta_{<\alpha,\,\alpha_0>}}
{2+\alpha^*_0(\alpha_0)\gamma_{\alpha_0}}\ \text{ for all non zero
}\alpha\text{ in }H(\tau)\vspace{0.3cm}\end{array}
\end{equation}

\section{Performing the computations}

\subsection{Multiplication by $p$ on an elliptic
curve}\label{computelliptic}

We refer the reader to \cite{Si} where it is explained how one can recover the group law on an elliptic curve $E$ via its geometry. More specifically, there are duplication and
addition formulae (see \cite[III,2]{Si}) given for an affine model
$y^2=x(x-1)(x-\mu)$ of $E$ as well as the division
polynomials in characteristic $p\geq 5$ (see \cite[Exercise
3.7]{Si}) which are much more convenient when implemented with a
computer. As the action of $\{\pm\}$ commutes with multiplication by
$p$, the latter induces a map $\P^1 \to\P^1$ defined by two homogeneous polynomials of degree $p$
(say $D$ and $N$) with $x \mapsto N(x^p)/D(x^p)$ and the map induced by Verschiebung is defined by $x \mapsto N(u)/D(u)$. The computations yield
\begin{lemma} In characteristic $p=3$,
\begin{equation*}N(u)=u(u+\mu(\mu+1))^2 \text{   and  }
D(u)=((\mu+1)u+\mu^2)^2\end{equation*} In characteristic
$p=5$,
\begin{equation*}
N(u)=
u\left[u^2-\mu(\mu+1)(\mu^2-\mu+1)u+\mu^4(\mu^2-\mu+1)\right]^2\end{equation*}
and
\begin{equation*}D(u)=\left[(\mu^2-\mu+1)\left[u^2-\mu^2(\mu+1)u\right]
+\mu^6\right]^2\end{equation*}
In characteristic $p=7$,
\begin{eqnarray*} N(u)& & =
u\left[u^3+2\mu(\mu+1)(\mu-2)(\mu-4)(\mu^2+3\mu+1)u^2\right.\\
& &
\hspace{0.7cm}\left.+\mu^4(\mu+1)^2\mu-2)(\mu-4)(\mu^2+1)u+
\mu^9(\mu+1)(\mu-2)(\mu-4)\right]^2\end{eqnarray*}
and
\begin{eqnarray*}D(u) =\left[(\mu+1)(\mu-2)(\mu-4)\left[u^3 \right.\right.& +& \mu^2(\mu+1)(\mu^2+1)u^2\\
& + &
\left.\left.\mu^6(\mu^2+3\mu+1)u\right]+\mu^{12}\right]^2\end{eqnarray*}
\end{lemma}
\begin{remark} \label{supersingellipcurve} These results are
consistent with the classification of supersingular elliptic curves in small characteristics (see \cite[Chapter IV, Examples 4.23.1, 4.23.2 and 4.23.3]{Ha}).\end{remark}

Choose non-zero scalars  $a$ and $b$
such that $\omega=-\Frac{a^4+b^4}{a^2b^2}$ is different from $\pm 2$ (in particular, $a\neq \pm b$).
There is a unique linear automorphism of $\P^1$ mapping $(a:\,b)$ to
0, $(a:\,-b)$ to $1$ and and $(b:\,a)$ to $\infty$. It maps
$(b:\,-a)$ to the point ($\mu:\,1)$ with
$
\mu=\left(\Frac{b^2+a^2}{2ab}\right)^2=\Frac{2-\omega}{4}$. Letting $\lambda_0,\,\lambda_1$ denote the corresponding pair of homogeneous
of $\P^1$, the homogeneous polynomials $Q_0$ and $Q_1$ corresponding to $N$ and $D$ under this linear transformation can be exchanged by the action of a suitable element of $E[2]$ and the computations yield expressions that depend only on $\omega$ and not on the choice of $a$ and $b$ as expected
\begin{lemma}\label{multpelliptic}
With the notations given above, one has\\
$\bullet \ p=3$. $
Q_0(\lambda_0,\,
\lambda_1)=\lambda_0^3-\omega\lambda_0\lambda_1^2$.\\
\noindent $\bullet \ p=5$. $Q_0=\lambda_0^5+\omega(\omega^2+2)\lambda_0^3\lambda_1^2+(\omega^2+2)
\lambda_0\lambda_1^4.$\\
\noindent $\bullet \ p=7$.
$Q_0=\lambda_0^7-2\omega(\omega^4-1)\lambda_0^5\lambda_1^2
+\omega^2(\omega^2-1)(\omega^2-2)
\lambda_0^3\lambda_1^4-\omega(\omega^2-1)\lambda_0\lambda_1^6
$.
\end{lemma}

\begin{remark} \label{consistencyg3} Using the expression of $\omega$ in terms of $\mu$,
the chart \eqref{formomegatau} as well as the remark
\ref{supersingellipcurve} allow one to give a more precise
description of the locus of the Kummer surfaces such that the
corresponding genus 2 curves only have ordinary Prym varieties.
Namely, viewing the set of Kummer surfaces as an open subset of the double covering of
the affine 3-space $\mathbb{A}^3$ given by equation
\eqref{relKummercoef}, one has to exclude the inverse image of a
finite set of affine planes in $\mathbb{A}^3$. Therefore, one can
easily say when a Coble quartic has associated Kummer surfaces such
that the corresponding genus 2 curves only have ordinary Prym
varieties (there is a finite set of linear relations that the
coefficients of the Coble quartic should not satisfy) and these form
a dense subset in the set of all Coble quartics. Together with the
Proposition \ref{Nakajima}, this ensures that a general genus 3
curve is ordinary, that all of its Prym varieties $P_\tau$ are
ordinary, and that if $Y_\tau$ is the genus 2 curve associated to
$P_\tau$, then all the Prym varieties of $Y_\tau$ also are ordinary.
\end{remark}

\subsection{Equations of $\widetilde{V}$ for $p=3$} Let us start with the genus 2 case and provide an alternative proof
to the following result of Laszlo and Pauly (\cite[Thm 6.1]{LP2} where the result is proven for any curve, not only a general one).

\begin{theorem}\label{LPg2p3} Let $X$ be a general smooth and projective curve of genus 2
over an algebraically closed filed of characteristic $3$.\\
(1) There is an embedding $\alpha : {\rm Kum}_{X} \hookrightarrow
|2\Theta_1|$ such that the equality of divisors in $|2\Theta_1|$
$$\widetilde{V}^{-1}({\rm Kum}_{X})= {\rm Kum}_{X_1}+2\alpha({\rm
Kum}_{X})$$ holds \emph{scheme-theoretically}.\\
(2) The cubic equations of $\widetilde{V}$ are given by the 4
partial derivatives of the quartic equation of the Kummer surface
$\alpha({\rm Kum}_X) \subseteq |2\Theta_1|$. In other words,
$\widetilde{V}$ is the polar map of the surface $\alpha({\rm
Kum}_X)$.\end{theorem} \begin{proof} Because $X$ is general, all its
Prym varieties are ordinary (Proposition \ref{Nakajima}). Let
$V_{0}$ be a generator for the $\hat{H}$-invariant part of
$\widetilde{V}^*(W)$ in ${\rm Sym}^p\,W_1$. For $g=2$ and $p=3$, it
comes under the form
\begin{equation*}V_0=y_{00}^3+a_{01}y_{00}y_{01}^2+a_{10}y_{00}y_{10}^2+a_{11}y_{00}y_{11}^2+b
y_{01}y_{10}y_{11}\end{equation*}  We use the Lemmas
\ref{reductionalpha0=0} and \ref{reductionalpha0neq0g2},
the chart \eqref{formomegatau} and the
Lemma \ref{multpelliptic} to obtain
\begin{equation}V_{0}=y_{00}^3+2k_{01}y_{00}y_{01}^2+2k_{10}y_{00}y_{10}^2+2k_{11}y_{00}y_{11}^2+2k_{00}
y_{01}y_{10}y_{11}\end{equation} Then, one can deduce the
$V_\alpha:=\widetilde{V}(x_\alpha)$ ($\alpha=01,\,10,\,11$) by
permuting suitably the coordinate functions $y_\bullet$ in $V_{0}$.
Notice that $V_{\alpha}$ is the partial (with respect to
$y_{\alpha}$) of a quartic surface with equation \begin{equation*}
S+ 2k_{00} P+k_{10}Q_{01}+k_{01}Q_{10}+k_{11}Q_{11}\end{equation*}
(with $S,\,P,\,Q_{01},\,Q_{10}$ and $Q_{11}$ as in Lemma \ref{Kummerdetails}) hence
isomorphic to ${\rm Kum}_X$. Thus, the second point above is proven.

 The inverse image
$\widetilde{V}^{-1}({\rm Kum}_{X})$ can be computed explicitly as it
is defined by the pull-back $\widetilde{V}^*(K_X)$ of the equation
\eqref{eqKummer}. In other words, a few more computations enable us
to recover the first assertion of the Theorem. Namely, one knows
(see Diagram \eqref{VS_XVJcommute}) that the equation $K_{X_1}$ of
${\rm Kum}_{X_1}$ divides $\widetilde{V}^*(K_X)$. The
exact quotient $\widetilde{V}^*(K_X)/K_{X_1}$ coincide with the square of
$K_X$.\end{proof}

This geometric interpretation of $\widetilde{V}$ in genus 2 and
characteristic 3 has an analogous interpretation of the unique
$\widetilde{V}$ lifting the generalized Verschiebung in genus 3 and
characteristic 3, namely the theorem \ref{LPgeneralizedg3}.

\begin{proof}[Proof of the Theorem \ref{LPgeneralizedg3}]
The proposition \ref{J[2]equivliftV} and the Remark
\ref{uniquenessg3p3} ensure the existence and the uniqueness of
$\widetilde{V} : |2\Theta_1| \dashrightarrow |2\Theta|$ lifting $V_2
: S_{X_1} \dashrightarrow S_X$ which is given by a system of 8
cubics. Because $X$ is general, the Remark \ref{consistencyg3} tells
us that the associated Prym varieties (that are ordinary by
Proposition \ref{Nakajima}) are general enough to have Kummer
surfaces for which the previous theorem is valid. Let $V_{0}$ be the
only $\hat{H}$-invariant cubic of this system of 8 cubics.
In the same way as in genus 2, we use the Lemmas
\ref{reductionalpha0=0} and \ref{reductionalpha0neq0g3}, the equations \eqref{CobletoKumalpha0=0} and \eqref{CobletoKumalpha0neq0} (instead of the chart \eqref{formomegatau}), and the Theorem \ref{LPg2p3} (instead of the
Lemma \ref{multpelliptic}) to obtain
\begin{eqnarray*} V_0 = y_{0}^3& +&
2\gamma_{001}y_{0}y_{001}^2+2\gamma_{010}y_{0}y_{010}^2
+2\gamma_{011}y_{0}y_{011}^2\\
&  +& 2\gamma_{100}y_{0}y_{100}^2+ 2\gamma_{101}y_{0}y_{101}^2
+2\gamma_{110}y_{0}y_{110}^2+2\gamma_{111}y_{0}y_{111}^2\\
&  +&
\delta_{001}y_{010}y_{100}y_{110}+\delta_{010}y_{001}y_{100}y_{101}
+\delta_{011}y_{011}y_{100}y_{111}\\
&  +& \delta_{100}y_{001}y_{010}y_{011}+
\delta_{101}y_{010}y_{101}y_{111}
+\delta_{110}y_{001}y_{110}y_{111}+\delta_{111}y_{011}y_{101}y_{110}
\end{eqnarray*}
Then, one deduces the $V_\alpha:=\widetilde{V}(x_\alpha)$ by
permuting suitably the coordinate functions $y_\bullet$ in $V_{0}$and checks that $V_{\alpha}$ is the partial (with respect to
$y_{\alpha}$) of a quartic hypersurface with equation
\begin{equation*}S^C +\sum
_{\alpha \in H-\{0\}} \gamma_\alpha Q^C_\alpha + \sum _{\alpha^* \in
\hat{H} -\{0\}} \delta_{\alpha^*} P^C_{\alpha^*}\end{equation*}
as in the Subsection \ref{sectCoblequartic},
hence isomorphic to ${\rm Cob}_X$. \end{proof}

\subsection{Equations of $\widetilde{V}$ for $p=5$ and $7$ and geometric consequences} For
these characteristics, we have only performed the calculations in the genus
2 case. The equations obtained are already huge and interpretation
requires (much time to) computational softwares. In much the same way as we proved the Theorem \ref{LPg2p3}, we prove the following

\begin{proposition} \label{eqVtildp5}Let $X$ be a general proper and smooth curve of genus 2
over an algebraically closed field of characteristic 5. There are
coordinate functions  $\{x_\alpha\}$ and $\{y_\alpha\}$ for
$|2\Theta|$ and $|2\Theta_1|$ respectively such that the Kummer
surface ${\rm Kum}_X$ in $|2\Theta|$ has an equation of the form
\eqref{eqKummer} and such that, if
$V_\alpha(\underline{y})=\widetilde{V}^*(x_\alpha)$, then
\begin{eqnarray*} V_{00}& = &
y_{00}^5+a_{1100}y_{00}^3y_{01}^2+a_{1010}y_{00}^3y_{10}^2+a_{1001}y_{00}^3y_{11}^2
+a_{0200}y_{00}y_{01}^4+a_{0110}y_{00}y_{01}^2y_{10}^2\\
& & +a_{0101}y_{00}y_{01}^2y_{11}^2+a_{0020}y_{00}y_{10}^4
+a_{0011}y_{00}y_{10}^2y_{11}^2+a_{0002}y_{00}y_{11}^4\\
& &
+b_{00}y_{00}^2y_{01}y_{10}y_{11}+b_{01}y_{01}^3y_{10}y_{11}+b_{10}y_{01}y_{10}^3y_{11}+b_{11}y_{01}y_{10}y_{11}^3
\end{eqnarray*}
with $$\begin{array}{lclcl}
a_{1100}=k_{01}(k_{01}^2+2), & & a_{1010}=k_{10}(k_{10}^2+2), & & a_{1001}=k_{11}(k_{11}^2+2),\\
a_{0200}=(k_{01}^2+2), & & a_{0020}=(k_{10}^2+2),& &
a_{0002}=(k_{11}^2+2),\\ \end{array}$$
$$
\begin{array}{rcl}
a_{0110}& = &3k_{11}(k_{00}^2+k_{11}^2)+k_{01}k_{10}(1-k_{11}^2),\\
a_{0101}& =& 3k_{10}(k_{00}^2+k_{10}^2)+k_{01}k_{11}(1-k_{10}^2),\\
a_{0011}& = & 3k_{01}(k_{00}^2+k_{01}^2)+k_{10}k_{11}(1-k_{01}^2),\\
b_{00}& = & 2k_{00}(k_{00}^2+1)-k_{00}k_{01}k_{10}k_{11},\\
\end{array}$$
$$b_{01}=k_{00}(k_{01}+3k_{10}k_{11}),\ b_{10}=k_{00}(k_{10}+3k_{01}k_{11}),\
b_{11}=k_{00}(k_{11}+3k_{01}k_{10})$$where the $k_i$ are the
coefficients of the equation \eqref{eqKummer} of ${\rm Kum}_X$. The
$V_\alpha$ ($\alpha=01,\,10,\,11$) can be deduced from $V_{00}$ by a
suitable permutation of the coordinate functions $y_i$, namely the
unique pairwise permutation that exchanges $y_{00}$ and $y_\alpha$.
\end{proposition}

\begin{remark} \label{eqVtildp7} In characteristic 7, one has the same kind of statement except that $V_0$ is now the sum of 30 monomials, whose coefficients are as above polynomials (of degree at most 7) in the parameters $k_\bullet$. \end{remark}

Let us recall some features of the geometry of the map
$\widetilde{V} : |2\Theta_1| \dashrightarrow |2\Theta|$ when $g=2$
and $p=3$ exhibited in \cite{LP2}. First, there is an
irreducible reduced hypersurface $H=\alpha({\rm Kum}_X)$ of degree
$2p-2=4$ in $S_{X_1}$ such that the equality of divisors
$$\widetilde{V}^{-1}({\rm Kum}_X)={\rm Kum}_{X_1} +2H$$ holds
in $S_{X_1}$ and such that base locus of $\widetilde{V}$ coincides
with the singular locus of $H$. Second, $H$ contains 16
curves (namely, conics containing 6 of the 16 singular points of
$H$) each of which was contracted by $\widetilde{V}$ on a singular
point of ${\rm Kum}_X$. Therefore, the inverse image of the singular
locus of ${\rm Kum}_X$ by $\widetilde{V}$ is 1-dimensional and
contains all the singular points of $H$.

We tried to check these properties in characteristic 5 and 7 in
exploiting the equations we could compute. Unfortunately,
most of these assertions require a Groebner basis computation which
seems beyond the capacities of the machines (or the software) used
when working with the generic curve (that is to say over the
extension
$$L:=\mathbb{F}_p(k_{01},\,k_{10},\,k_{11})[k_{00}]/(k_{00}^2-k_{01}^2-k_{10}^2-k_{11}^2+k_{01}k_{10}k_{11}+4)$$
of the prime field $\mathbb{F}_p$ as field of coefficients).
Still, we could prove Proposition
\ref{LPgeneralizedg2}.

\begin{proof}[Proof of the proposition \ref{LPgeneralizedg2}] Define the field
$L$ as above and $R$ as the $L$-graded algebra generated by
$y_{00},\,y_{01},\,y_{10}$ and $y_{11}$. The homogeneous polynomials
$V_{00},\,V_{01},\,V_{10}$ and $V_{11}$ define a endomorphism of
$L$-graded algebras $\widetilde{V}^* : R \to R$. Letting $K$ (resp.
$K_1$) be the equation \eqref{eqKummer} of the Kummer surface ${\rm
Kum}_X$ in $|2\Theta|$ (resp. ${\rm Kum}_{X_1}$ in $|2\Theta_1|$,
which is obtained by raising the coefficients of \eqref{eqKummer} to
the power $p$), one checks that $K_1$ divides $\widetilde{V}^*(K)$.
Letting $Q$ be the exact quotient $\widetilde{V}^*(K)/K_1$, one
checks that it is a square $S^2$ and that $S$ is
irreducible.\end{proof}

Although the computation does not end for the generic curve, we
checked the other assertions for a one hundred plus particular
curves in each characteristic ($p=5$ and 7), randomly choosing the
coefficients $k_{01},\,k_{10},\,k_{11}$ in $\mathbb{F}_{p^{10}}$.
The following pattern arises, with no exception once put aside the
(non-relevant) cases contradicting the inequalities
\eqref{coefconstraints}. The base locus of
$\widetilde{V}$ is contained, scheme-theoretically, in the singular
locus of $H$, which is 0-dimensional and has total length 96 (resp.
304) in characteristic 5 (resp. 7). This singular locus is itself in
the stable locus of $S_{X_1}$ as well as in the inverse image of the
singular points of ${\rm Kum}_X$ which again defines a 1-dimensional
subset of $H$. The inclusion of the base locus of $\widetilde{V}$ in
the singular locus of $H$ is strict and, unfortunately, its
reducedness is too expensive to be checked out by the computation.\\

Of course, one would like to find a geometric (and characteristic
free) proof of these facts that we believe are true for any general
curve and any odd characteristic.

\end{document}